\documentclass[12pt,article]{amsart}

\usepackage{amsmath,amssymb,epsfig,color, amsthm, mathrsfs}
\usepackage{wrapfig,lipsum,floatrow,sidecap,comment}
\usepackage{graphicx}
\usepackage{multirow}
\usepackage{hhline}
\usepackage{url}
\usepackage{soul}
\graphicspath{ {images/} }
\allowdisplaybreaks
\def\sideremark#1{\ifvmode\leavevmode\fi\vadjust{\vbox
to0pt{\vss \hbox to 0pt{\hskip\hsize\hskip1em
\vbox{\hsize2cm\tiny\raggedright\pretolerance10000
\noindent#1\hfill}\hss}\vbox to8pt{\vfil}\vss}}}

\theoremstyle{plain} 
\newcommand{\Z}{{\mathbb Z}}
\newcommand{\C}{{\mathbb C}}

\newcommand{\T}{{\mathbb T}}

\renewcommand{\L}{{\mathcal L}}
\newcommand{\MP}{{\mathfrak P}}

\def\p{\prime}
\def\pp{{\prime\prime}}
\def\xxi{{\bf{\xi}}}
\def\e{{\bf{e}}}

\def\KN{\K\setminus \{\0\}}
\def\LN{\L\setminus \{\0\}}
\def\C{\mathbb{C}}

\def\x{{\bf{x}}}
\def\y{{\bf{y}}}
\def\z{{\bf{z}}}

\def\k{{\bf{k}}}
\def\l{{\bf{l}}}
\def\n{{\bf{n}}}
\def\m{{\bf{m}}}

\def\0{{\bf{0}}}
\def\e{{\bf{e}}}
\def\Z{\mathbb{Z}}
\def\T{\mathbb{T}}

\def\O{\Omega}
\def\f{{\bf{f}}}
\def\F{\mathcal{F}}
\def\K{\mathcal{K}}
\def\L{\mathcal{L}}

\newtheorem {theo} {\bf Theorem} [section]

\newtheorem {prop} [theo] {\bf Proposition}

\newtheorem {lem} [theo] {\bf Lemma}

\newtheorem {rem} [theo] {\bf Remark}
\newcommand{\R}{\mathbb R\,}

\def\F{\mathcal{F}}

\numberwithin{equation}{section}
\begin{document}
\title{Gabor single-frame and multi-frame  multipliers in any given dimension}
\author{Yuanan Diao$^\dagger$}
\address{$^\ddagger$ Department of Mathematics and Statistics\\
University of North Carolina Charlotte, Charlotte, NC 28223, USA}
\email{ydiao@uncc.edu}
\author{Deguang Han$^\ddagger$}
\address{$^\dagger$ Department of Mathematics\\
University of Central Florida\\ Orlando, FL,
32816, USA}
\email{deguang.han@ucf.edu}
\author{Zhongyan Li$^\sharp$}
\address{$^\sharp$ Department of Mathematics and Physics\\
North China Electric Power University, Beijing, 102206, China}
\email{lzhongy@ncepu.edu.cn}

\date{\today}
\thanks{ Deguang Han acknowledges the support from NSF under the grant DMS-1712602; Zhongyan Li acknowledges the support from the
National Natural Science Foundation of China (Grant No. 11571107).}
\maketitle
\begin{abstract}
Functional  Gabor single-frame or multi-frame generator multipliers are the matrices of  function entries  that preserve Parseval Gabor single-frame or multi-frame generators. An interesting and natural  question is how to  characterize all such multipliers. This question has been answered for several special cases including the  case of single-frame generators in two  dimensions and the case of multi-frame generators in one-dimension. In  this paper we  completely characterize multipliers  for  Gabor single-frame and multi-frame generators with respect to  separable time-frequency lattices in any given dimension. Our approach is general and applies to the previously known cases as well.
\end{abstract}

\section{Introduction}

This paper is a continuation of a project on characterizing the frame generator multipliers under various settings, a topic that was initially motivated by the work   of Dai and Larson \cite{DL} on wandering vector multipliers and the WUTAM paper on basic properties of wavelets \cite{WC}.  Representative publications  resulted from this project include \cite{G Han, Han-TAMS, LDDX, LHan, LHan-LAA}. Let  $A$ and $B$ be two nonsingular real matrices such that $A\Z^d$ and $B\Z^d$ are both full-rank lattices in $\R^d$.
We say that $G(\x) =(g_{1}(\x), ... , g_{\gamma}(\x) )^{\tau}$,  with $\tau$ being the transpose and $g_j(\x)\in L^2(\R^d)$ for each $j$,  is a {\em Parseval (or normalized tight) Gabor multi-frame generator} of length $\gamma$ for $L^{2}(\R^d)$ (for the separable  time-frequency lattice $A\Z^{d}\times B\Z^{d}$) if $\{e^{2\pi i\langle B\m,\x\rangle}g_j(\x-A\n):\,\m,\n\in\Z^d, 1\le j\le \gamma\}$ is a normalized tight frame, {\em i.e.},
\begin{equation}\label{Parseval_eq}
\sum_{1\le j\le \gamma}\sum_{\m,\n\in\Z^d}|\langle f, {e^{2\pi i\langle B\m,\x\rangle}}g_j(\x-A\n)\rangle|^2=\|f\|^2
\end{equation}
for all $f(\x)\in L^2(\R^d)$. In the special case that $\gamma=1$, $G(\x)=g_1(\x)$ is also called a {\em Parseval Gabor single-frame generator} or simply a {\em Parseval Gabor frame generator}. However for the sake of convenience in this paper a Parseval Gabor single-frame generator will simply be regarded as a Parseval Gabor multi-frame generator with length $\gamma=1$. Gabor multi-frames in higher dimensions play important roles in many applications (\cite{FH1,FH2,JF,LS,ZY,ZM}).

\medskip
A functional  matrix
$M(\x)=(f_{ij}(\x))_{\gamma\times \gamma}$ with $f_{ij}(\x) \in L^{\infty}(\R^d)$ is called a functional (matrix) Gabor multi-frame multiplier if $H(\x) = M(\x)G(\x)$ is a Parseval Gabor multi-frame generator for $L^2(\R^d)$ whenever $G=(g_1,g_2,\cdots,g_\gamma)^{\tau}$ is. Functional (matrix) Gabor multi-frame multipliers provide a useful tool in the study of Parseval Gabor multi-frames. As such, it is an interesting and important question to ask how they can be characterized. This question has been answered for the special cases of $d=1$ with any $\gamma\ge 1$ (\cite{G Han, LHan-17}), and $d=2$ with $\gamma=1$ \cite{LiD}. In this paper, we provide a complete characterization for functional (matrix) Gabor multi-frame multipliers with any number $\gamma\ge 1$ of generators at any dimension $d\ge 1$. Our result shall contain all the previously obtained results in \cite{G Han,LiD, LHan-17} with a unified approach. More specifically, we have proven the following theorem.

\medskip
\begin{theo}\label{main-thm} Let $M(\x)=(f_{ij}(\x))_{\gamma\times \gamma}$ ($\gamma\ge 1$) with $f_{ij}(\x)\in L^{\infty}(\R^d)$. Then $M(\x)$  is a functional matrix Gabor multi-frame multiplier for the time-frequency lattice $A\Z^{d}\times B\Z^{d}$ if and only if the following three conditions are satisfied:

\noindent\smallskip
(1) $M(\x)$ is unitary for a.e. $\x\in\R^d;$

\noindent\smallskip
(2) For any $\n\in\Z^d\setminus\{\0\}$, $M^*(\x)M(\x+(B^{\tau})^{-1}\n)$ equals $\lambda_{\n}(\x)I$ ({\em a.e.} $\x\in\R^d$) for  some unimodular scalar-valued function (that depends only on $\n$) $\lambda_{\n}(\x)$, where $I$ is the identity matrix and $M^*$ denotes the conjugate transpose of $M$.

\noindent\smallskip
(3) $\lambda_{\n}(\x)$ is $A\Z^d$-periodic (as a function of $\x$).
\end{theo}

\medskip
Notice that in the case of $\gamma=1$, a functional matrix Gabor multi-frame multiplier is a scalar function $h(\x)$ hence Theorem \ref{main-thm} has a simpler form:
\begin{theo}\label{main-thm2} A scalar function $h(\x)\in L^{\infty}(\R^d)$ is a functional  Gabor multiplier for the time-frequency lattice $A\Z^{d}\times B\Z^{d}$ if and only if the following two conditions hold:

\noindent\smallskip
(1) $h(\x)$ is unimodular for a.e. $\x\in\R^d;$

\noindent\smallskip
(2) For any $\n\in\Z^d\setminus\{\0\}$, $h(\x)\overline{h(\x+(B^{\tau})^{-1}\n)}$ is $A\Z^d$-periodic.
\end{theo}

\medskip
We should note that it is known that a Gabor multi-frame generator of length $\gamma$ exists if and only if $|\det(AB)|\le \gamma<|\det(AB)|+1$ (\cite{JD}). Thus $|\det(AB)|\leq 1$ corresponds to the case of $\gamma=1$ which is the (implied) condition for Theorem \ref{main-thm2} and $1< |\det(AB)|$ corresponds to the case of $\gamma=\lceil |\det(AB)|\rceil \ge 2$ which is the (implied) condition for Theorem \ref{main-thm}.

\medskip
The following well known characterization for normalized tight Gabor multi-frame generators \cite{Casa} plays a key role in our proofs of Theorem \ref{main-thm} and Theorem \ref{main-thm2}.

\medskip
\begin{prop}\label{Lem1} \cite{Casa}
Let $A,\,B$ be nonsingular matrices with $|\det(A)|=a,\,|\det(B)|=b$, and $g_1, g_2, \cdots, g_\gamma\in L^2(\R^d).$ Then $G = (g_1, ... , g_{\gamma})^{\tau}$  is a Parseval Gabor multi-frame generator for $L^2(\R^d)$ if and only if the following identities hold
(for a.e. $\x\in\R^d$):
\begin{eqnarray}
\label{cond1}&&\sum_{\n\in\Z^d}\langle G(\x-A\n), G(\x-A\n)\rangle=b;\\
\label{cond2}&&\sum_{\n\in\Z^d}\langle G(\x-A\n), G(\x+(B^{\tau})^{-1}\l-A\n)\rangle=0, \,\forall \ \l\in\Z^d\setminus\{\0\}.
\end{eqnarray}
\end{prop}

\medskip
Using Proposition \ref{Lem1}, the sufficient part of Theorem \ref{main-thm} can be easily proven. Indeed, let $M(\x)$ be a functional matrix satisfying the conditions (1)--(3) in Theorem \ref{main-thm}, and $G(\x)=(g_1,\cdots,g_\gamma)^{\tau}$ be an arbitrary Parseval Gabor multi-frame generator. Denote $H(\x)=M(\x)G(\x)=(\eta_1(\x),\cdots,\eta_\gamma(\x))^{\tau}$. Since $M(\x)$ is unitary for any  $\x\in\R^d$ {\em a.e.}, it is obvious that
$$\sum_{\n\in\Z^d}\langle
H(\x-A\n),H(\x-A\n)\rangle=\sum_{\n\in\Z^d}\langle G(\x-\n),G(\x-\n)\rangle,
$$
hence (\ref{cond1}) holds.
Furthermore, by conditions (2) and (3) we have
\begin{eqnarray*}
&&\sum_{\n\in\Z^d}\langle H(\x-A\n),H(\x-A\n+(B^{\tau})^{-1}\l)\rangle\\
&=&\sum_{\n\in\Z^d}\langle G(\x-A\n), M^*(\x-A\n)M(\x-A\n+(B^{\tau})^{-1}\l)G(\x-A\n+(B^{\tau})^{-1}\l)\rangle\\
&=&\sum_{\n\in\Z^d}\langle G(\x-A\n), M^*(\x)M(\x+(B^{\tau})^{-1}\l)G(\x-A\n+(B^{\tau})^{-1}\l)\rangle\\
&=&\overline{\lambda_\l(\x)}\sum_{\n\in\Z^d}\langle G(\x-A\n), G(\x-A\n+(B^{\tau})^{-1}\l)\rangle=0
\end{eqnarray*}
for any $\x\in\R^d$ {\em a.e.} and any $\l\neq \0$. Hence $H(\x)=M(\x)G(\x)$ is a Parseval Gabor multi-frame generator. Thus the rest of the paper is devoted to the proof of the necessary part of Theorem \ref{main-thm}.

\medskip
The rest of the paper is organized as follows. In the next section, we show that Theorem \ref{main-thm} can be proven under a simplified setting. In Section 3 we provide some necessary background knowledge regarding the lattice tiling and packing of $\R^d$. In Section 4 we prove Theorem \ref{main-thm} for the case of $\gamma=1$ with any $d\ge 1$. In Section 5 we prove Theorem \ref{main-thm} for the case of $\gamma>1$ with any $d\ge 1$.

\section{Auxiliary Simplifications}

Let  $A$ and $B$ be nonsingular real valued $d\times d$ matrices, and $\gamma$ be the integer satisfying $|\det(AB)|\le \gamma<|\det(AB)|+1$. Let $P$, $Q$ be any two $d\times d$ matrices with integer entries and $|\det(P)|=|\det(Q)|=1$. (Note that we will be making specific choices later for $P$, $Q$ later depending on our needs, but the statement here holds for any such $P$, $Q$.) Denote the matrix $(PB^\tau AQ)^{-1}$ by $D$. This implies that  $AQD=(B^\tau)^{-1}P^{-1}$. For any set of functions $\{\tilde{g}_1(\x),...,\tilde{g}_\gamma(\x)\}$ such that $\tilde{g}_j\in L^2(\R^d)$, define $g_j(\x)$ by $g_j(\x)=\tilde{g}_j((AQ)^{-1}\x)=\tilde{g}_j(\z)$ where $\z=(AQ)^{-1}\x$.

\medskip
\begin{lem}\label{Lemma2}
The following two statements hold:\\
\noindent
(1) $\tilde{G}(\x)=(\tilde{g}_1(\x),...,\tilde{g}_\gamma(\x))^\tau$ is a Parseval Gabor multi-frame generator for the time-frequency lattice $\Z^d\times (D^\tau)^{-1}\Z^d$ if and only if $G(\z)=(g_1(\z),...,g_\gamma(\z))^\tau$ is a Parseval Gabor multi-frame generator for the time-frequency lattice $A\Z^d\times B\Z^d$.\\
\noindent
(2) Let $\tilde{M}(\z)$ be a $\gamma\times \gamma$ functional matrix multiplier and define $M(\x)=\tilde{M}(\z)$ with $\z=(AQ)^{-1}\x$. If $\tilde{M}^*(\z){\tilde{M}(\z+D\k)}$ is $\Z^d$-periodic for any $\z\in \R^d$ and $\k\in\Z^d\setminus\{\0\}$, then $M^*(\x){M(\x+(B^{\tau})^{-1}\k)}$ is $A\Z^d$-periodic for any $\x\in \R^d$ and $\k\in \Z^d\setminus\{\0\}$. Moreover, if $\tilde{M}^*(\z)\tilde{M}(\z+D\k)=\tilde{\lambda}_{\k}(\z)I$ for some scalar function $\tilde{\lambda}_{\k}(\z)$, then $M^*(\x){M(\x+(B^{\tau})^{-1}\k^\p)}=\lambda_{\k^\p}(\x)I$ with $\lambda_{\k^\p}(\x)=\tilde{\lambda}_{\k}(\z)$.
\end{lem}

\begin{proof} (1) We leave it to our reader to verify that the system $\{e^{2\pi i\langle (D^{\tau})^{-1}\k,\z\rangle}\tilde{g}_j(\z-\l): \l,\k\in \Z^d, 1\le j\le \gamma\}$ can be rewritten as $\{e^{2\pi i\langle B\k,\x\rangle}g_j(\x-A\l): \l,\k\in\Z^d, 1\le j\le \gamma\}$. The statement then follows easily by verifying (\ref{Parseval_eq}) for the latter system using the variable substitution $\x=(AQ)\z$ in the integrals.

\medskip\noindent
(2) For any $\l_1$, $\k_1(\neq 0)\in\Z^d$, let $\l=Q^{-1}\l_1\in \L$ and $\k=P\k_1\in \Z^d$, then for any $\x\in \R^d$ we have
\begin{eqnarray*}
&&M^*(\x-A\l_1)M(\x-A\l_1+(B^{\tau})^{-1}\k_1))\\
&=&
M^*(AQ(\z-\l))M(AQ(\z-\l+D\k))=\tilde{M}^*(\z-\l)\tilde{M}(\z-\l+D\k)\\
&=&
\tilde{M}^*(\z)\tilde{M}(\z+D\k)=M^*(AQ(\z))M(AQ(\z+D\k))
=
M^*(\x)M(\x+(B^{\tau})^{-1}\k_1).
\end{eqnarray*}
This shows that $M^*(\x)M(\x+(B^{\tau})^{-1}\k_1))$ is $A\Z^d$-periodic for any $\x\in \R^d$. The last statement in (2) is obvious.
\end{proof}

\medskip
Lemma \ref{Lemma2} implies that in order to prove Theorem \ref{main-thm}, we only need to consider a special case of it, namely when $A=I_{d\times d}$ and $(B^T)^{-1}=D$, where $D$ is of the form $(PB^TAQ)^{-1}$ for any $P$, $Q$ with integer entries and $|\det(P)|=|\det(Q)|=1$. Under this setting,
 it is necessary that $0<|\det(D^{-1})|=d_0\le \gamma <d_0+1$ and equations (\ref{cond1}) and (\ref{cond2}) become
\begin{eqnarray}
\sum_{\l\in\L}\langle G(\x-\l),G(\x-\l)\rangle&=&d_0, \label{cond12}\\
\sum_{\l\in\L}\langle G(\x-\l),G(\x-\l-\k)\rangle&=&0,\ \forall\ \k\in \K\setminus \{\0\}. \label{cond22}
\end{eqnarray}

\section{The tiling and packing of $\R^d$ by $\L$ and $\K$}\label{Tiling}

Let us introduce a few key concepts  first. Let $\F$  be any full rank lattice of $\R^d$.  A measurable set $E$ is said to {\em pack} $\R^d$ by $\F$ if $E\cap (E+\f)=\emptyset $ for any nontrivial $\f\in \F$. If $E$ packs $\R^d$ by $\F$ and also satisfies the condition $\R^d=\cup_{\f\in \F}(E+\f)$,  then we say that $E$ {\em tiles} $\R^d$ by $\F$. In this case $E$ is called a {\it tile} or a {\it fundamental domain} of $\F$. For two measurable sets $S_1$ and $S_2$ that pack $\R^d$ by $\F$, we say that $S_1$ and $S_2$ are {\em $\F$-equivalent} if $\cup_{\f\in \F}(S_1+\f) = \cup_{\f\in \F}(S_2+\f)$, and we say that $S_1$ and $S_2$ are {\em $\F$-disjoint} if $S_1\cap (S_2+\f)=\emptyset$ for any nontrivial $\f\in \F$.

\medskip
The materials in this section heavily rely on the work \cite{Han W}, more specifically the proofs in the  sequence of lemmas there that lead to the proof of \cite[Theorem 1.2]{Han W}, which states that if $|\det(D)|\le 1$, then there exists a measurable set that tiles $\R^d$ by $\K=D\Z^d$ and packs $\R^d$ by $\L=\Z^d$. From this point on, the lattices $\L$, $\K$ always mean $\Z^d$, $D\Z^d$ respectively unless otherwise noted. The following long remark summarizes the results (with slight modifications) extracted from \cite{Han W} that are necessary  for us to prove Theorem \ref{main-thm}.

\medskip
\begin{rem}\label{long_remark}{\em
Consider the group $\T^d=\R^d/\Z^d$ with $\O=[0,1)^d$ a  representative  set of the group. Let $\pi: \ \R^d\longrightarrow \T^d$ be the projection map and consider $\pi((B^\tau A)^{-1}\Z^d)$. $\overline{\pi((B^\tau A)^{-1}\Z^d)}$ is a closed subgroup of $\T^d$, hence $\overline{\pi((B^\tau A)^{-1}\Z^d)}=S\oplus F$ for some rational
subspace $S$ and finite set $F$ \cite[Lemma 2.1]{Han W}. The proof of \cite[Theorem 1.2]{Han W} is divided into three cases: \textbf{Case 1}: $S=\T^d$; \textbf{Case 2}: $S=\{\emptyset\}$; and \textbf{Case 3}: $S\not=\T^d$ and $S\not=\{\emptyset\}$. We will follow these cases to make the  choices for $P$, $Q$ and to extract the information  we need. Let $\k_j=D\l_j$ where $\{\l_1,\l_2,...,\l_d\}$ is the standard basis for $\L$. Let $\gamma$ be the unique integer satisfying $d_0\le \gamma<d_0+1$.

\medskip\noindent
{\bf Case 1.} In this case we can simply choose $P=Q=I_{d\times d}$ in Lemma \ref{Lemma2}. There are two sub cases to consider here: (i) $d_0=1/|\det(D)|=|\det(AB)|$ is rational and (ii) $d_0$ is irrational.

\medskip\noindent
(i) $d_0=p/q$ with $(p,q)=1$. In this case we can partition $D\O$ into $M_2$ parallelepipeds of the same volume $\mu_0$ where $M_2$ can be any integer multiple of $q$. We have $M_2\mu_0=q/p=1/d_0$ hence $\mu_0=q/(M_2p)=1/(N p)$ where $N q=M_2$. On the other hand, we can partition $\O$ into $M_1=N p$ rectangles such that each rectangle also has volume $\mu_0$.  Denote these partitions of $\O$ and $D\O$ by $\MP$ and $\MP^\p$ respectively, and arbitrarily order and name the ones in $\MP$ as $C_1$, $C_2$, ..., $C_{M_1}$, and the ones in $\MP^\p$ as $C^\p_1$, $C^\p_2$, ..., $C^\p_{M_2}$. \cite[Corollary 2.3]{Han W} assures that for any pair of $C_i$ and $C_j^\p$, there exists a measurable set $J(C_i,C^\p_j)$ that is $\L$-equivalent to $C_i$ and $\K$-equivalent to $C_j^\p$ (we say $J(C_i,C^\p_j)$ is a {\em matching} of $C_i$  and $C_{j}^\p$). In particular, if there exists a rectangle $C$ such that $C\subset C_i$ and $-\l_0+C\subset -\k_0+C_j^\p$ for some $\l_0\in \L$ and $\k_0\in \K$, then $-\l_0+C$ can be selected as part of $J(C_i,C^\p_j)$.  Recall that $M_1=p N$ and $M_2=qN$. If $q=1$, then $d_0=p=\gamma$ is an integer and $M_1=\gamma M_2$. This means in this case we can divide the rectangles in $\MP$ into $\gamma$ groups $F_1$, $F_2$, ..., $F_\gamma$ such that each group contains $M_2$ rectangles. If $p>q>1$ (this happens when $d_0>1$), then $p=(\gamma -1)q+r$ for some positive integer $r<q$ (otherwise $p=\gamma q$ contradicts the condition that $(p,q)=1$). It follows that  $M_1=Np=(\gamma -1)Nq+Nr=(\gamma-1)M_2+Nr$. Thus in this case we can divide the rectangles in $\MP$ into $\gamma$ groups $F_1$, $F_2$, ..., $F_\gamma$ such that each group $F_j$ with $j\ge 2$ contains $M_2$ rectangles, and $F_1$ contains the remaining $Nr<Nq=M_2$ rectangles. Finally, if $q>p\ge 1$ (that is, $d_0<1$), then we have $M_1<M_2$.

\medskip\noindent
(ii) $d_0$ is irrational. Here we need to consider the cases $d_0>1$ and $d_0<1$ separately.

\medskip
First consider the case $d_0>1$. We have $d_0=\gamma -\delta$ for some positive constant $\delta<1$.  In this case we can still partition $D\O$ into $M_2$ parallelepipeds (denoted by $C_j^\p$'s) of the same volume $\mu_0$ where $M_2$ can be any arbitrarily chosen large positive integer, in particular, we will choose it  large enough so that  $(1-\delta)M_2>1$. This time it is not possible to partition $\O$ into rectangles such that each rectangle also has the same volume $\mu_0$ since $\mu_0=1/(M_2d_0)$ is irrational, however this can be done if we allow one of these rectangles to have volume less than $\mu_0$. Thus if $M_1$ is the total number of rectangles in  $\MP$ named and ordered as $C_i$'s as before, we can assume that all $C_i$'s have volume $\mu_0$ except that $C_{M_1}$ has volume $\mu^\p$ which is less than $\mu_0$. We leave it to our reader to verify that in this case $M_1-(\gamma -1)M_2=(1-\delta)M_2+1-(\mu^\p/\mu_0)>(1-\delta)M_2>1$ by the choice  of $M_2$. This means that we can again divide the rectangles in $\MP$ into $\gamma$ groups $F_1$, $F_2$, ..., $F_\gamma$ such that each group $F_j$ with $j\ge 2$ contains $M_2$ rectangles, and $F_1$ contains the remaining rectangles including $C_{M_1}$ which has volume $\mu^\p$. By the above inequality we see that $F_1$ contain at least two rectangles, hence it also contains at least one rectangle that has volume $\mu_0$. The statement in (i) about  $J(C_i,C^\p_j)$ applies if $i\not=M_1$. For $i=M_1$, $C_{M_1}$ can be matched to any parallelepiped of volume $\mu^\p$ that is a subset of any of the $C_j^\p$'s.

\medskip
Now consider the case $d_0<1$. In this case we first partition $\O$ into $M_1$ rectangular parallelepipeds of the same volume $\mu_0$ where $M_1$ can be arbitrarily large. For example we can partition $\O$ into $N^d$ small cubes of side length $1/N$ where $N>0$ can be any arbitrarily chosen integer, obtaining $M_1=N^d$ cubes, each with volume $\mu_0=1/N^d$. On the other hand, we can partition $D\O$ into $M_2=N^d$ parallelepipeds such that each parallelepiped has volume $1/(d_0N^d)>\mu_0$. Denote these partitions of $\O$ and $D\O$ by $\MP$ and $\MP^\p$ respectively, and arbitrarily order and name the ones in $\MP$ as $C_1$, $C_2$, ..., $C_{M_1}$, and the ones in $\MP^\p$ as $C^\p_1$, $C^\p_2$, ..., $C^\p_{M_1}$. \cite[Corollary 2.3]{Han W} assures that for any pair of $C_i$ and $C_j^\p$ ($1\le i,j\le M_1$), there exists a measurable set $J(C_i,C^\p_j)$ that is $\L$-equivalent to $C_i$ and $\K$-equivalent to a subset of $C_j^\p$. We call $J(C_i,C^\p_j)$ a {\em matching} of $C_i$  and $C_{j}^\p$. In particular, if there exists a rectangular parallelepiped $C$ such that $C\subset C_i$ and $-\l_0+C\subset -\k_0+C_j^\p$ for some $\l_0\in \L$ and $\k_0\in \K$, then $-\l_0+C$ can be selected as part of $J(C_i,C^\p_j)$.

\medskip\noindent
{\bf Case 2.} In this case $P$ and $Q$ can be chosen so that $D$ is a diagonal matrix with positive rational entries
$\{p_1/q_1,p_2/q_2,...,p_d/q_d\}$ (with ${\rm gcd}(p_i,q_i)=1$). Here we need to choose $\MP$ and $\MP^\p$ specifically so that the elements in them are rectangular parallelepipeds whose sides are parallel to $\l_1$, ..., $\l_d$ with corresponding side lengths of $1/q_1$, $1/q_2$, ..., $1/q_d$. Again, name and order (arbitrarily) the rectangles in $\O$ and $D\O$ by $C_i$'s and $C_j^\p$'s with $1\le i\le M_1=q_1q_2\cdots q_d$ and $1\le j\le M_2=p_1p_2\cdots p_d$. Then for any pair $C_i$ and $C_j^\p$, there exists a rectangle   $J(C_i,C_j^\p)$ such that $J(C_i,C_j^\p)=\l_0+C_i=\k_0+C$ for some $\l_0\in \L$ and $\k_0\in \K$.  In particular, if $C_i$ is paired with $C_j^\p$ and $-\l_0+C_i=-\k_0+ C_j^\p$ for some $\l_0\in \L$ and  $\k_0\in K$, then $J(C_i,C_j^\p)=-\l+C_i$. Since $d_0$ is rational, by a discussion similar to Case 1 (i), if $d_0>1$ then the rectangles in $\MP$ can be divided into $\gamma$ groups $F_1$, $F_2$, ..., $F_\gamma$ such that each group $F_j$ with $j\ge 2$ contains $M_2$ rectangles, and $F_1$ contains at least one and at most $M_2$ rectangles. On the other hand, if $d_0\le 1$, then $M_1\le M_2$.

\medskip\noindent
{\bf Case 3.} In this case $P$ and $Q$ can be chosen so that
$D=\begin{bmatrix}D_1&B_0\\0&D_2\end{bmatrix}$, where $D_2$ is a $d_2\times d_2$ diagonal matrix with positive rational entries $\{p_1/q_1,...,p_{d_2}/q_{d_2}\}$  (${\rm gcd}(p_i,q_i)=1$ for each $i$)  and $\begin{bmatrix}D_1& B_0\end{bmatrix}\Z^d$ (mod 1) is dense in $[0,1)^{d_1}$, $d_1+d_2=d$.
For the sake of  convenience  let $\{\l_1^\p,\l_2^\p,...,\l^\p_{d_1}\}$  be the standard basis  of $\R^{d_1}$, $\K_1$ be the lattice spanned  by  $\{\k_1^\p,\k_2^\p,...,\k^\p_{d_1}\}$ where $\k^\p_j=D_1\l_j^\p$ ($1\le j\le d_1$). Similarly, let $\{\l_1^\pp,\l_2^\pp,...,\l^\pp_{d_2}\}$  be the standard basis  of $\R^{d_2}$, $\K_2$ be the lattice spanned  by  $\{\k_1^\pp,\k_2^\pp,...,\k^\pp_{d_2}\}$ where $\k^\pp_j=D_2\l_j^\pp=(p_j/q_j)\l_j^\pp$ ($1\le j\le d_2$). Let $\O_1=[0,1)^{d_1}$, $\O_2=[0,1)^{d_2}$  so that $\O=\O_1\times \O_2$. It is apparent that $D_1\O_1$ is a fundamental domain of $\K_1$ and $D_2\O_2$ is a fundamental domain of $\K_2$. It is less apparent  that $\tilde{\O}=(D_1\O_1)\times (D_2\O_2)$ is also a fundamental domain of $\K$ \cite{Han W}. So instead of using $D\O$, we will use $\tilde{\O}=(D_1\O_1)\times (D_2\O_2)$ instead. The advantage of this is that it allows us to obtain a partition of $\tilde{\O}$ by partitioning $D_1\O_1$  and $D_2\O_2$ separately as described below.

\medskip\noindent
(i) $|\det(D_1)|$ is rational. Similar to Case 1 (i), we can partition $D_1\O_1$ into parallelepipeds $\Delta_j^\p$'s ($1\le j\le M^\p_2$) of the same volume $\mu_0$ and partition $\O_1$ into rectangular parallelepipeds $\Delta_i$'s ($1\le i\le M^\p_1$) with volume $\mu_0$ ($\mu_0$ can be chosen to be as small as we want). On the other hand, we partition $\O_2$ and $D_2\O_2$ into rectangles $R_i$'s ($1\le i\le q_1  q_2\cdots q_{d_2}$) and $R_j^\p$'s ($1\le j\le p_1  p_2\cdots p_{d_2}$) whose sides are parallel to $\l^\pp_1$, ..., $\l^\pp_{d_2}$ with corresponding side lengths of $1/q_1$, $1/q_2$, ..., $1/q_{d_2}$. Then the set $\{\Delta_j^\p\times R_i^\p:\ 1\le j\le M_2^\p, 1\le i\le p_1 p_2\cdots p_{d_2}\}$ is a partition $\MP^\p$ of $\tilde{\O}$ whose total number of elements is $M_2=M_2^\p p_1 p_2\cdots p_{d_2}$, and the set $\{\Delta_j\times R_i:\ 1\le j\le M_1^\p, 1\le i\le q_1 q_2\cdots q_{d_2}\}$ is a partition $\MP$ of ${\O}$ whose total number of elements is $M_1=M_1^\p q_1 q_2\cdots q_{d_2}$. By a slightly modified version of \cite[Sub Lemma 5]{Han W}, for any $C_{ij}=\Delta_i\times R_j\in \MP$ and any $C^\p_{i^\p j^\p}=\Delta_{i^\p}^\p\times R_{j^\p}^\p\in \MP^\p$, there exists a measurable set $J(C_{ij},C^\p_{i^\p j^\p})$ that is $\L$-equivalent to $C_{ij}$ and $\K$-equivalent to $C^\p_{i^\p j^\p}$. In particular, if $\Delta_i$ contains a small rectangle $\Delta_0$ such that $-\l_0+\Delta_0\times R_j \subset -\k_0+\Delta_{i^\p}^\p\times R_{j^\p}^\p=-\k_0+C^\p_{i^\p j^\p}$, then $-\l_0+\Delta_0\times R_j$ can be selected as part of $J(C_{ij},C^\p_{i^\p j^\p})$. In  the case that $d_0>1$, the discussion in Case 1 (i) applies to $M_1-(\gamma-1)M_2$ here (since $d_0$ is rational). That is, we can divide the elements in $\MP$ into $\gamma$ groups $F_1$, $F_2$, ..., $F_\gamma$ such that each group $F_j$ with $j\ge 2$ contains $M_2$ elements, and $F_1$ contains at least one and at most $M_2$ elements. On the other hand, if $d_0\le 1$, then $M_1\le M_2$.

\medskip\noindent
(ii) $|\det(D_1)|$ is irrational hence $d_0$ is also irrational. Similar to Case 1(ii), we also need to consider the cases $d_0>1$ and $d_0<1$ separately.

\medskip
If $d_0>1$, let $\delta=\gamma-d_0>0$. Similar to (i) above, we partition $\O_2$ and $D_2\O_2$ into rectangles $R_i$'s ($1\le i\le q_1  q_2\cdots q_{d_2}$) and $R_j^\p$'s ($1\le j\le p_1  p_2\cdots p_{d_2}$) whose sides are parallel to the $\l^\pp_1$, ..., $\l^\pp_{d_2}$ coordinates with corresponding side lengths of $1/q_1$, $1/q_2$, ..., $1/q_{d_2}$. Similar to Case 1 (ii), we can partition $D_1\O_1$ into parallelepipeds $\Delta_j^\p$'s ($1\le j\le M^\p_2$) of the same volume $\mu^\p_0$ and partition $\O_1$ into rectangles $\Delta_i$'s ($1\le i\le M^\p_1$) with volume $\mu^\p_0$ ($\mu_0^\p$ can be chosen to be as small as we want since we can choose $M_2^\p$ as large as we want), with the exception that $0<\mu^\pp=\mu(\Delta_{M_1^\p})<\mu_0^\p$. Similar to (i) above, the set $\{\Delta_j^\p\times R_i^\p:\ 1\le j\le M_2^\p, 1\le i\le p_1 p_2\cdots p_{d_2}\}$ is a partition $\MP^\p$ of $\tilde{\O}$ whose total number of elements is $M_2=M_2^\p p_1 p_2\cdots p_{d_2}$, and the set $\{\Delta_j\times R_i:\ 1\le j\le M_1^\p, 1\le i\le q_1 q_2\cdots q_{d_2}\}$ is a partition $\MP$ of ${\O}$ whose total number of elements is $M_1=M_1^\p q_1 q_2\cdots q_{d_2}$. The difference here is that all the elements in these partitions have measure $\mu_0=\mu_0^\p/(q_1 q_2\cdots q_{d_2})$, except that the elements $\Delta_{M_1^\p}\times R_i$ (for any $1\le i\le q_1 q_2\cdots q_{d_2}$) have measure $\mu^\pp/(q_1 q_2\cdots q_{d_2})$. However, the inequality $M_1-(\gamma -1)M_2>(1-\delta)M_2$ still holds in this case as one can check, hence we will have $(1-\delta)M_2>q_1 q_2\cdots q_{d_2}+1$ if $M_2^\p$ is large enough. This ensures that we can place all the elements $\Delta_{M_1^\p}\times R_i$ into the group $F_1$ as described in Case 1 (ii), and $F_1$ will still contain at least one element whose measure is $\mu_0$. The statement in (i) about the matching of two elements with measure $\mu_0$ applies. A set of the form  $\Delta_{M_1^\p}\times R_i$, on the other hand, can be matched with a set $\Delta\times R_j^\p$, where $\Delta$ is a properly chosen parallelepiped contained in a $\Delta_i^\p$ with a measure $\mu^\pp$.

\medskip
In the case that $d_0<1$, we first partition $\O_2$ into rectangular parallelepipeds whose sides are parallel to $\l^\pp_1$, ..., $\l^\pp_{d_2}$ with corresponding side lengths of $1/q_1$, $1/q_2$, ..., $1/q_{d_2}$, and partition $\O_1$ into $M_1^\p=N^{d_1}p_1p_2\cdots p_{d_2}$ (congruent) rectangular parallelepipeds whose sides are parallel to $\l^\p_1$, ..., $\l^\p_{d_1}$ with corresponding side lengths of $1/N$, $1/N$, ..., $1/N$, $1/(Np_1p_2\cdots p_{d_2})$, where $N$ can be any arbitrarily chosen positive integer. Combining these partitions yields a partition $\MP$ of $\O=\O_1\times \O_2$. If we name and number the rectangular parallelepipeds in the partition of $\O_1$ as $\Delta_1$, $\Delta_2$, ..., $\Delta_{M_1^\p}$, name and number the rectangular parallelepipeds in the partition of $\O_2$ as $R_1$, $R_2$, ..., $R_{M_1^\pp}$ where $M_1^\pp=q_1 q_2\cdots q_{d_2}$, then  $\MP=\{C_{ij}=\Delta_i\times R_j:\ 1\le i\le M_1^\p, 1\le j\le M_1^\pp\}$, and $M_1=M_1^\p M_1^\pp=N^{d_1}(p_1p_2\cdots p_{d_2})(q_1 q_2\cdots q_{d_2})$ is the total number of partition elements in $\MP$. Notice that the volume of each $\Delta_j$ is $\mu_0^\p=1/M_1^\p=1/(N^{d_1}p_1p_2\cdots p_{d_2})$ and the volume of each $R_i$ is $\mu_0^\pp=1/(q_1 q_2\cdots q_{d_2})$.

\medskip
Next, we partition $D_1\O_1$ into $M_2^\p=N^{d_1}q_1 q_2\cdots q_{d_2}$ congruent parallelepipeds. The volume of each of these  parallelepipeds is $\mu(D_1\O_1)/M_2^\p=|\det(D_1)|/M_2^\p\ge 1/(N^{d_1} p_1p_2\cdots p_{d_2})=\mu_0^\p$. We name and order these parallelepipeds as
$\Delta_i^\p$'s ($1\le i\le M^\p_2$). We now partition $D_2\O_2$ into rectangular parallelepipeds that are congruent to the $R_i$'s and name/order them as $R_j^\p$'s ($1\le j\le p_1  p_2\cdots p_{d_2}$). Combining these two partitions yields a partition $\MP^\p$ of $\tilde{\O}=D\O_1\times D\O_2$: $\MP^\p=\{C^\p_{i^\p j^\p}=\Delta_{i^\p}^\p\times R_{j^\p}^\p:\ 1\le {i^\p}\le M_2^\p, 1\le {j^\p}\le M_2^\pp\}$ where $M_2^\pp=p_1 p_2\cdots p_{d_2}$. Notice that the total number of partition elements in $\MP^\p$ is $M_2=M_2^\p M_2^\pp=N^{d_1}(p_1p_2\cdots p_{d_2})(q_1 q_2\cdots q_{d_2})=M_1$ but $\mu(\Delta_{i^\p}^\p\times R_{j^\p}^\p)\ge \mu_0^\p/(q_1 q_2\cdots q_{d_2})=\mu(\Delta_i\times R_j)$.
Again by a slightly modified version of \cite[Sub Lemma 5]{Han W}, for any $C_{ij}=\Delta_i\times R_j\in \MP$ and any $C^\p_{i^\p j^\p}=\Delta_{i^\p}^\p\times R_{j^\p}^\p\in \MP^\p$, there exists a measurable set $J(C_{ij},C^\p_{i^\p j^\p})$ that is $\L$-equivalent to $C_{ij}$ and $\K$-equivalent to a subset of $C^\p_{i^\p j^\p}$. In particular, if $\Delta_i$ contains a small rectangular parallelepiped $\Delta_0$ such that $-\l_0+\Delta_0\times R_j \subset -\k_0+\Delta_{i^\p}^\p\times R_{j^\p}^\p=-\k_0+C^\p_{i^\p j^\p}$, then $-\l_0+\Delta_0\times R_j$ can be selected as part of $J(C_{ij},C^\p_{i^\p j^\p})$.

\medskip
This ends Remark \ref{long_remark}.
}
\end{rem}

\medskip
Notice  that in Remark \ref{long_remark}, we used $D\O$ as the  fundamental domain of $\K$ in  Cases 1 and 2,  but used $\tilde{\O}$  as the  fundamental domain of $\K$ in  Case 3. For the sake of simplicity, let us denote  them by $\O^\p$ with the understanding that it means either $D\O$ or $\tilde{\O}$ depending on which case $\K$ belongs to. Let us call the pair of partitions $\MP$, $\MP^\p$ of $\O$ and $\O^\p$ with the properties discussed in Remark \ref{long_remark} {\em nice partition pair}.

\medskip
\begin{rem}\label{translation}{\em
We note that in the above discussion, $\O^\p$ can be replaced by any $\K$-translation of $\O^\p$. Similarly, $\O$ can be replaced by any $\L$-translation of $\O$. Thus without loss of generality, one can always assume that $\O\cap \O^\p$ is non-empty and contain interior points.
}
\end{rem}

\section{The proof of Theorem  \ref{main-thm}, Part I}

\medskip
We now proceed to prove Theorem  \ref{main-thm2}, namely the special case $\gamma=1$ of Theorem  \ref{main-thm}. Here $d_0=|\det(D^{-1})|\le 1$ so $|\det(D)|=1/d_0\ge 1$. In  this case, the  functional multiplier is a scalar function $h(\x)$ and we need to prove that $h(\x)\overline{h(\x-\k)}=h(\x-\l)\overline{h(\x-\l-\k)}$ for any $\x\in \R^d$ {\em a.e.}, and for any $\l\in \L$ and $\k\in \KN$.

\medskip
\subsection{A few useful lemmas} First, for the very special case of $D=I$, we have the following lemma, which follows by a generalized version of the  proof of \cite[4.3 Type III Case]{LiD}.

\begin{lem}\label{D=I}
If $D=I$ is the identity matrix, then Theorem \ref{main-thm} holds.
\end{lem}

\medskip
We also have the following lemma,  whose proof can be found in \cite{LiD}.

\medskip
\begin{lem}\label{Lemma6}\cite[Lemma 2.4]{LiD}
Let $h$ be a functional Gabor frame multiplier and let $\l\in \L$, $\k\in \KN$ be any given pair of vectors,  then $h(\x)\overline{h(\x- \k)}=h(\x- \l)\overline{h(\x- \l-\k)}$ for any $\x\in \R^d$ if one of the following conditions holds:\\
\noindent
(i) There exist disjoint and measurable sets $E_1$, $E_2$, $E_3$, $E_4$ and $E_5$, with $E_2$ being a rectangular parallelepiped, such that $E_3=-\k+E_2$, $E_4=-\l+E_2$, $E_5=-\l-\k+E_2$, and $E_1\cup E_2\cup E_3$ tiles $\R^d$ by $\L$ while $E_1\cup E_2\cup E_4$ packs $\R^d$ by $\K$;
\\
\noindent
(ii) $\k\in \L$ and there exist disjoint and measurable sets $E_1$, $E_2$, $E_3$, $E_4$ and $E_5$, with $E_2$ being a rectangular parallelepiped, such that $E_3=-\k+E_2$, $E_4=-\l+E_2$, $E_5=-\l-\k+E_2$, $E_1\cup E_2$ tiles $\R^d$ by $\L$ and $E_1\cup E_2\cup E_4$ packs $\R^d$ by $\K$.
\end{lem}

\medskip
Notice that the given condition in (i) above implies that $\k\not\in \L$ and $\l\not\in \K$, while the given condition in (ii) implies that $\l\not\in \K$.

\medskip
\begin{lem}\label{Lemma7}
Let $h$ be a functional Gabor frame multiplier and let $\l,\  \l^\p\in \L$, $\k,\ \k^\p\in \KN$, then the following statements hold:\\
\noindent
(i) If the following two equations hold for any $\x\in \R^d$:
\begin{eqnarray*}
h(\x)\overline{h(\x-\k)}&=&h(\x- \l^\p)\overline{h(\x- \l^\p-\k)},\\
h(\x)\overline{h(\x-\k)}&=&h(\x+\l^\p-\l)\overline{h(\x+\l^\p-\l-\k)},
\end{eqnarray*}
 then 
$
h(\x)\overline{h(\x-\k)}=h(\x- \l)\overline{h(\x- \l-\k)}
$
for any $\x\in \R^d$.
\\
\noindent
(ii) If the following two equations hold  for any $\x\in \R^d$: 
\begin{eqnarray*}
h(\x)\overline{h(\x-\k^\p)}&=&h(\x- \l)\overline{h(\x- \l-\k^\p)},\\
h(\x)\overline{h(\x-\k^\p-\k)}&=&h(\x-\l)\overline{h(\x-\l-\k^\p-\k)},
\end{eqnarray*}
 then 
 $
 h(\x)\overline{h(\x-\k)}=h(\x- \l)\overline{h(\x- \l-\k)}
 $
 for any $\x\in \R^d$.
\end{lem}

\begin{proof}
(i) Substituting $\x+\l^\p$ by $\x$ (by a slight abuse of notation) on both sides of the equation
$$
h(\x)\overline{h(\x-\k)}=h(\x+\l^\p-\l)\overline{h(\x+\l^\p-\l-\k)},
$$
we obtain
$$
 h(\x-\l^\p)\overline{h(\x-\l^\p-\k)}=h(\x-\l)\overline{h(\x-\l-\k)}.
$$
But the left side of the above equation is $h(\x)\overline{h(\x-\k)}$ by the given condition.

\medskip
(ii) Since $h$  is unimodular, the given condition leads to
$$
h(\x)\overline{h(\x- \l)}=h(\x- \k^\p)\overline{h(\x- \l- \k^\p)}
$$
and
$$
h(\x)\overline{h(\x- \l)}=h(\x- \k^\p-\k)\overline{h(\x- \l- \k^\p-\k)}.
$$
Thus we have
$$
h(\x- \k^\p)\overline{h(\x- \l- \k^\p)}=h(\x- \k^\p-\k)\overline{h(\x- \l- \k^\p-\k)}.
$$
Substituting $\x-\k^\p$ by $\x$ on both sides of the above equation (again by a slight abuse of notation), we obtain
$$
h(\x)\overline{h(\x- \l)}=h(\x-\k)\overline{h(\x- \l-\k)}
$$
for any $\x\in \R^d$ {\em a.e.}
\end{proof}

\medskip
\begin{lem}\label{choice_lemma}
For any given $\k$ and $\l$, if $\k\not\in \L$ and $\l\not\in \K$, then there exist disjoint and measurable sets $E_1$, $E_2$, $E_3$, $E_4$ and $E_5$, with $E_2$ being a rectangular parallelepiped, such that $E_3=-\k+E_2$, $E_4=-\l+E_2$, $E_5=-\l-\k+E_2$, and $E_1\cup E_2\cup E_3$ tiles $\R^d$ by $\L$ while $E_1\cup E_2\cup E_4$ packs $\R^d$ by $\K$. It follows that
$h(\x)\overline{h(\x- \l)}=h(\x-\k)\overline{h(\x- \l-\k)}$ by Lemma \ref{Lemma6}.
\end{lem}

\medskip
\begin{proof} We will prove the given statement by discussing the three different cases of $D$ given in Remark \ref{long_remark}. The general strategy in each case is to choose a suitable nice partition pair $\MP$, $\MP^\p$ of $\O$ and  $\O^\p$ as described in Remark \ref{long_remark} with the following additional properties: (a) there exists a rectangular parallelepiped (which is chosen as our set $E_2$) such that $E_2\subset C_0\cap C_0^\p$ for some $C_0\in \MP$ and $C^\p_0\in \MP^\p$; (b) $-\k-\l_0+E_2\subset C_1\in \MP$, for some $\l_0\in \L$ and $C_1\not=C_0$; (c) $-\l-\k_0+E_2\subset C^\p_1\in \MP^\p$, for some $\k_0\in \K$ and $C^\p_1\not=C^\p_0$. We then assign a one to one correspondence between the elements of $\MP$ and $\MP^\p$ in an arbitrary way to be used to define the matchings, except that  $C_0$ is matched with $C_0^\p$ with  the choice of $E_2\subset J(C_0,C_0^\p)$ and $C_1$ is matched with  $C_1^\p$ with the choice that $\k-\l_0+C$ is matched with $\l-\k_0+C$ in $J(C_1,C_1^\p)$. Denote the union of the matchings by $S$. $S$ tiles $\R^d$ by $\L$ and packs $\R^d$ by $\K$ by our construction. Let the portion of $J(C_1,C_1^\p)$ that is $\L$-equivalent to $-\k-\l_0+C$ and $\K$-equivalent to $-\l-\k_0+C$ by $J_0$, and define $E_1=S\setminus (E_2\cup J_0)$, $E_3=-\k+E_2$, $E_4=-\l+E_2$ and $E_5=-\l-\k+E_2$. Since  $-\l_0+E_3=-\k-\l_0+E_2$, $E_3$ is $\L$-equivalent to $J_0$ hence $E_1\cup E_2\cup E_3$ is $\L$-equivalent to $S$. Thus $E_1\cup E_2\cup E_3$ tiles $\R^d$ by $\L$. Similarly, $E_4$ is $\K$-equivalent to $J_0$ hence $E_1\cup E_2\cup E_4$ is $\K$-equivalent to $S$, so $E_1\cup E_2\cup E_4$ packs $\R^d$ by $\K$. Thus in the following cases we only need to show that a nice partition pair $\MP$, $\MP^\p$ of $\O$ and  $\O^\p$ that satisfies conditions (a) to (c) exists.

\medskip\noindent
Case 1. Let $\x_0$ be an interior point of $\O\cap D\O$ near the origin. Since $\O$ is a fundamental domain of $\L$, there exists $\l_0\in \L$ such that $\x_0-\k-\l_0\not=\x_0$ is also an interior point of $\O$ (since $-\k\not\in \L$). Similarly, there exists $\k_0\in \K$ such that $\x_0-\l-\k_0\not=\x_0$ is an interior point of $D\O$. It follows that we can choose a nice partition pair $\MP$, $\MP^\p$ of $\O$ and $D\O$ such that the diameters of the polytopes in them are smaller than $\min\{|-\k-\l_0|,|-\l-\k_0|\}$. WLOG we can assume that there exist $C_{i_0}\in \MP$ and $C_{j_0}^\p\in \MP^\p$ such that $\x_0$ is an interior point of $C_{i_0}\cap C_{j_0}^\p$, $\x_0-\k-\l_0$ is an interior point of another $C_{i_1}$ in $\MP$, and $\x_0-\l-\k_0$ is an interior point of another $C^\p_{j_1}$  in $\MP^\p$ (otherwise we can replace $\x_0$ by a suitably chosen point very close to it). It follows that there exists a small rectangular parallelepiped $E_2$ such that $E_2\subset C_{i_0}\cap C_{j_0}^\p$, $\k-\l_0+C\subset C_{i_1}$ and $\l-\k_0+C\subset C^\p_{j_1}$.

\medskip\noindent
Case 2. In this case $\MP$ and $\MP^\p$ are as described in Case 2 of Remark \ref{long_remark} with the rectangular parallelepipeds named and ordered in an arbitrary way. Let $E_2$ be the common rectangular parallelepiped of $\MP$ and $\MP^\p$ containing the origin (say $E_2=C_{i_0}=C_{j_0}^\p$). Similar to the proof of Case 1, there exist $\l_0\in \L$ and $\k_0\in \K$ such that $-\k-\l_0+E_2=C_{i_1}\in \MP$ for some $C_{i_1}\in \MP$, and  $-\l-\k_0+E_2=C^\p_{j_1}$ for some $C^\p_{j_1}\in \MP^\p$. $i_1\not=i_0$ and $j_1\not=j_0$ since $-\l\not\in \K$ and $-\k\not\in\L$.

\medskip\noindent
Case 3. Recall that in this case $D=\begin{bmatrix}D_1&B_0\\0&D_2\end{bmatrix}$, where $D_2$ is a $d_2\times d_2$ diagonal matrix with positive rational entries, $\begin{bmatrix}D_1& B_0\end{bmatrix}\Z^d$ (mod 1) is dense in $[0,1)^{d_1}$, and $\tilde{\O}=D_1\O_1\times D_2\O_2$ is a fundamental domain of $\K$. There exist $\l_0\in \L$ and $\k_0\in \K$ such that $-\k-\l_0=\begin{bmatrix}\y_1\\ \y_2\end{bmatrix}$ with $\y_1\in \O_1$, $\y_2\in \O_2$ and at least one of them is not $\0$, $-\l-\k_0=\begin{bmatrix}\z_1\\ \z_2\end{bmatrix}$ with $\z_1\in D_1\O_1$, $\z_2\in D_2\O_2$ and at least one of them is not $\0$. Thus the nice partition pair $\MP$, $\MP^\p$ of $\O$ and $\tilde{O}$ obtained by combining a suitably chosen nice partition pair $\MP_1$, $\MP_1^\p$ of $\O_1$ and $D_1\O_1$, and the nice partition pair $\MP_2$, $\MP_2^\p$ of $\O_2$ and $D_2\O_2$ as described in Case 2 of Remark \ref{long_remark}, we will have the following: there exists a small rectangular parallelepiped $C\subset \O_1\cap D_1\O_1$ such that $E_2=C\times R$, where $R=R_{i_0^\p}=R_{j_0^\p}^\p$ is the common rectangular parallelepiped of $\MP_2$ and $\MP^\p_2$ containing the origin, satisfies the condition that $E_2\subset C_{i_0,i_0^\p}\cap C^\p_{j_0,j_0^\p}$, $C_{i_0,i_0^\p}=\Delta_{i_0}\times R_{i_0^\p}\in \MP$ and $C^\p_{j_0,j_0^\p}=\Delta_{j_0}^\p\times R_{j_0^\p}^\p\in \MP^\p$, $-\k-\l_0+E_2\subset C_{i_1,i_1^\p}=\Delta_{i_1}\times R_{i_1^\p}\in \MP$  ($C_{i_1,i_1^\p}\not=C_{i_0,i_0^\p}$), $-\l-\k_0+E_2\subset C^\p_{j_1,j_1^\p}=\Delta^\p_{j_1}\times R^\p_{j_1^\p}\in \MP^\p$ ($C^\p_{j_1,j_1^\p}\not=C^\p_{j_0,j_0^\p}$).
\end{proof}

\medskip
\subsection{The main proof} We shall consider the following three cases: I. $\L\not\subset \K$ and $\K\not\subset \L$; II. $\L\subset \K$; III. $\K\subset \L$ but $\L\not\subset \K$.

\bigskip\noindent
\textbf{I. $\L\not\subset \K$ and $\K\not\subset \L$} Let $\l\in \L$ and $\k\in\KN$ be  any two vectors. If $\l\not\in\K$ and $\k\not\in \L$, then by Lemma \ref{choice_lemma} we have $h(\x)\overline{h(\x- \l)}=h(\x-\k)\overline{h(\x- \l-\k)}$. If $\l\not\in \L$ and $\k\in\KN$, then choose any $\k^\p\in\KN$ such that $\k^\p\not\in \L$, and apply Lemma \ref{choice_lemma} to the pairs $(\l, \k^\p)$ and $(\l, \k^\p+\k)$ (note that $\k^\p+\k\not\in\L$). Similarly, if $\l\in \L$ and $\k\not\in\KN$, then choose any $\l^\p\in\L$ such that $\l^\p\not\in \K$, and apply Lemma \ref{choice_lemma}  to the pairs $(\l^\p,\k)$ and $(-\l^\p+\l,\k)$ (note that $-\l^\p+\l\not\in\K$). Finally, if $\l\in \K$ and $\k\in\L$, then choose any $\l^\p\in\L$ such that $\l^\p\not\in \K$ and $\k^\p\in\KN$ such that $\k^\p\not\in \L$. The application of Lemma \ref{choice_lemma} to the pairs $(\l^\p,\k^\p)$ and $(-\l^\p+\l, \k^\p)$ leads to $h(\x)\overline{h(\x- \l)}=h(\x-\k^\p)\overline{h(\x- \l-\k^\p)}$, and to the pairs $(\l^\p,\k^\p+\k)$ and $(-\l^\p+\l,\k^\p+\k)$ leads to $h(\x)\overline{h(\x- \l)}=h(\x-\k^\p-\k)\overline{h(\x- \l-\k^\p-\k)}$. The result now follows by Lemma \ref{Lemma7} (ii).

\medskip\noindent
\textbf{II. $\L\subset \K$} Notice that $\begin{bmatrix}\l_1&\l_2 & \cdots & \l_d\end{bmatrix}=I$ is the identity  matrix, so we have $I=DU$ for some $d\times d$ matrix $U$ with integer entries since each $\l_j$ is a  linear combination of $\k_1=D\l_1$, $\k_2=D\l_2$, ..., $\k_d=D\l_d$ with integer coefficients. It follows that  $D^{-1}=U$ is a matrix with integer entries hence $|\det(D^{-1})|\ge 1$ since $|\det(D^{-1})|$ is positive and is an integer. However we also have $|\det(D^{-1})|=d_0\le 1$. Thus $|\det(D^{-1})|= 1$ and it follows that $D$ itself is a matrix with integer entries and $|\det(D)|=1$. By Case 2 of Remark \ref{long_remark}, we must have $D=I$ and the result follows from Lemma \ref{D=I}.

\medskip\noindent
\textbf{III. $\K\subset \L$ but $\L\not\subset \K$} Since $\L\not\subset \K$, there exists $\l^\p\in \L$ such that $\l^\p\not\in \K$. In this case $\O$ tiles $\R^d$ by $\L$ and packs $\R^d$ by $\K$, and $E_1=\emptyset$, $E_2=\O$, $E_3=-\k+E_2$, $E_4=-\l^\p+E_2$, $E_5=-\l^\p-\k+E_2$ satisfy condition (ii) of Lemma \ref{Lemma6} since $-\l^\p+E_2$ is $\K$-disjoint from $E_2$. It follows that
$$
h(\x)\overline{h(\x-\k)}=h(\x- \l^\p)\overline{h(\x- \l^\p-\k)}
$$
for any $\x\in \R^d$ {\em a.e.}
Similarly, the sets $E^\p_1=\emptyset$, $E^\p_2=\O$, $E^\p_3=-\k+E^\p_2$, $E^\p_4=-\l+\l^\p+E^\p_2$, $E^\p_5=-\l+\l^\p-\k+E^\p_2$ also satisfy condition (ii) of Lemma \ref{Lemma6} since $-\l+\l^\p+E_2$ is also $\K$-disjoint from $E_2$. Thus we have
$$
h(\x)\overline{h(\x-\k)}=h(\x+ \l^\p-\l)\overline{h(\x+\l^\p-\l-\k)}.
$$
The result now follows from Lemma \ref{Lemma7} (i).

\section{The proof of Theorem  \ref{main-thm}, Part 2}\label{Sec_Proof}

We now prove the remaining case of Theorem \ref{main-thm}, namely the case $\gamma>1$, which is uniquely determined by $1<d_0\le \gamma<d_0+1$.

\medskip
\subsection{A few more lemmas}

\medskip
\begin{lem}\label{Lem3}
Let $D$ be as given in Remark \ref{long_remark},  $\gamma$ be the unique integer determined by $d_0\le \gamma <d_0+1$ where $d_0=1/|\det(D)|>1$. Then there exist measurable sets $E_1$, $E_2$, ..., $E_\gamma$ such that \\
\noindent
(1) $\mu(E_j)=|\det(D)|=1/d_0$ for $2\le j\le \gamma$ and $0<\mu(E_1)\le 1/d_0$;\\
\noindent
(2) $E_i$ and $E_j$ are $\L$-disjoint if $i\not=j$ and $\cup_{1\le j\le \gamma}E_j$ tiles $\R^d$ by $\L$;\\
\noindent
(3) $E_j$ tiles $\R^d$ by $\K$ for each $j\ge 2$ and $E_1$ packs $\R^d$ by $\K$.
\end{lem}

\begin{proof}
Let $\MP$ and $\MP^\p$ be as defined in Remark \ref{long_remark} and $F_1$, $F_2$, ..., $F_\gamma$ be the groups defined there. Assign an arbitrary one to one correspondence between the elements of $F_j$ and the elements of $\MP^\p$ ($j\ge 2$) so they are paired according to this one to one correspondence. The union of the matchings of these pairs is a measurable set $E_j$ that is $\L$-equivalent to the union of  the sets in $F_j$ and $\K$-equivalent to $\O^\p$. For $F_1$ we assign an arbitrary one to one mapping from $F_1$ to $\MP^\p$, and each pair leads to a matching of an element in $F_1$ and a subset of its corresponding element in $\MP^\p$. The union of these matchings is a measurable set $E_1$ that is $\L$-equivalent to the  union of the sets in $F_1$ and $\K$-equivalent to a subset of $\O^\p$. $E_i$ and  $E_j$  are $\L$-disjoint since $F_i$ and $F_j$ contain disjoint subsets of $\O$ if $i\not=j$, and $\cup_{1\le j\le \gamma}E_j$ is $\L$-equivalent to $\O$ hence tiles $\R^d$ by $\L$.
\end{proof}

\medskip
\begin{lem}\label{choice_lemma2}
For any given $\k\in \KN$, the following statements hold.

\medskip
(1) If $\k\not\in \L$, then there exists a  choice of $E_1$, ..., $E_\gamma$  such that $E_1$ contains a  rectangular parallelepiped $C$ such that $-\k+C\subset E_2$;

\medskip
(2) If $\k\in \L$, then there exists a choice of $E_1$, ..., $E_\gamma$ such that $E_1$ contains a rectangular parallelepiped $C$ such that $C^\p=-\k_{i_0}+C\subset E_2$ for some $\k_{i_0}=D\l_{i_0}$.
\end{lem}

\medskip
\begin{proof} (1) We will prove this by discussing the three different cases of $D$ given in Remark \ref{long_remark}.

\medskip\noindent
Case 1. Let $\l_0\in \L$ be such that $-\k-\l_0=\pi(-\k)\in \O$. Since $\k\not\in \L$, $\pi(-\k)\not=\0$. Let $\x_0$ be an interior point of $\O$ such that $\x_0+\pi(-\k)$ is also an interior point of $\O$. Choose a nice partition pair $\MP$, $\MP^\p$ of $\O$ and $D\O$ such that the diameters of the polytopes in them are smaller than $|\pi(-\k)|$. WLOG we can assume that $\x_0$ is an interior point of $D\O$ as well by Remark \ref{translation}. Furthermore, we can assume that there exist $C_{i_0}\in \MP$ and $C_{j_0}^\p\in \MP^\p$ such that $\x_0$ is an interior point of $C_{i_0}\cap C_{j_0}^\p$, and $\x_0+\pi(-\k)$ is also an interior point of another $C_{i_1}$  (otherwise we can replace $\x_0$ by a point very close to it). It follows that there exists a small rectangle $C$ such that $C\subset C_{i_0}\cap C_{j_0}^\p$ and $\pi(-\k)+C\subset C_{i_1}$. We now place $C_{i_0}$ in the group $F_1$ and $C_{i_1}$ in the group $F_2$ ($F_1$ and $F_2$ are as defined in Remark \ref{long_remark}), and match both $C_{i_0}$ and $C_{i_1}$ to $C_{j_0}^\p$. By Remark \ref{long_remark}, we can choose $C$ to be a subset of the matching of $C_{i_0}$ and $C_{j_0}^\p$ since it is both $\L$-equivalent and $\K$-equivalent to itself. Similarly, we can choose $-\k+C$ to be a subset of the matching of  $C_{i_1}$ and $C_{j_0}^\p$, since it is $\K$-equivalent to $C\subset C_{j_0}^\p$ and is $\L$-equivalent to  $\pi(-\k)+C\subset C_{i_1}$ (since $-\k+C=\l_0+(\pi(-\k)+C)\subset \l_0+C_{i_1}$). Thus $C\subset E_1$ and $-\k+C\subset E_2$.

\medskip\noindent
Case 2. Let $\l_0\in \L$ be such that $-\k-\l_0=\pi(-\k)\in \O$. Let $C^\p$ be the common rectangle of $\MP$ and $\MP^\p$ containing the origin. In this case $\pi(-\k)+C^\p$ is another rectangle of $\MP$. Similar to the discussion in Case 1, we can choose to have $C\subset E_1$ and $-\k+C\subset E_2$.

\medskip\noindent
Case 3. Let $\l_0\in \L$ be such that $-\k-\l_0=\pi(-\k)\in \O$. Similar to the discussion of Case 1, if we choose the diameters of the parallelepipeds in $\MP^\p$ small enough, then there exist $C_{ij}=\Delta_i\times R_j\in \MP$ and $C^\p_{i^\p j^\p}=\Delta_{i^\p}^\p\times R_{j^\p}^\p\in \MP^\p$ such that $R_j=R_{j^\p}^\p$, $\Delta_i\cap \Delta_{i^\p}^\p$ contains a small rectangle $C_0$, and $\pi(-\k)+C_0\times R_j\subset C_{i_1j_1}=\Delta_{i_1}\times R_{j_1}\in \MP$. We place $C_{ij}$ in the group $F_1$ and $C_{i_1j_1}$ in the group $F_2$, and match both $C_{ij}$ and $C_{i_1j_1}$ to $C^\p_{i^\p j^\p}$. Again by Remark \ref{long_remark}, we can choose $C=C_0\times R_i$ to be a subset of the matching of $C_{ij}$ and $C^\p_{i^\p j^\p}$ since it is both $\L$-equivalent and $\K$-equivalent to itself. Similarly, we can choose $-\k+C$ to be a subset of the matching of  $C_{i_1j_1}$ and $C^\p_{i^\p j^\p}$, since it is $\K$-equivalent to $C\subset C^\p_{i^\p j^\p}$ and is $\L$-equivalent to  $\pi(-\k)+C\subset C_{i_1j_1}$. Thus $C\subset E_1$ and $-\k+C\subset E_2$.

\medskip
(2) Notice that at least one of the $\k_j$'s has length less than $1$. Let $\k_{j_0}$ be one such. The  discussions in  the  above apply here with $\pi(-\k)$ replaced by $\k_{j_0}$. It follows that in each case there exists a small rectangle $C$ such that $C\subset E_1$ and $\k_{j_0}+C\subset E_2$.
\end{proof}

\medskip
\subsection{The main proof.}
We now proceed to prove Theorem \ref{main-thm} for the case of $\gamma\ge 2$, under the simplified setting of $A=I_{d\times d}$ and $(B^T)^{-1}=D$, where $D$ is of the form $(PB^TAQ)^{-1}$ as described in Remark \ref{long_remark}.
The proof is divided into three parts. In Part 1 we show that if $M(\x)$ is a Parseval Gabor multi-frame multiplier, then $M(\x)$ is unitary. In Part 2 we show that $M^*(\x)M^*(\x-\k)=\lambda_\k(\x)I$ for any $\k\in \K\setminus \{\0\}$, where $\lambda_\k(\x)$ is  a scalar function that depends only on $\k$ and $\x$. In the last part we show that $M^*(\x)M^*(\x-\k)$ is  $\L$-periodic.

\begin{rem}\label{E_j_remark}{\em
Notice that the discussions in the last section can be applied to any translation of $\O$ (together with the set $\O^\p$, of course). Thus in order to verify that equations (\ref{cond12}) and (\ref{cond22}) hold for any $\x\in \R^d$ {\em a.e.}, we only need to verify them for any $\x\in \O$ {\em a.e.} We should stress that the statement here is different from the statement of Remark \ref{translation}.
}
\end{rem}

\medskip\noindent
\textbf{Part 1.} Let $M(\x)$ be a $\gamma \times\gamma$ functional matrix Gabor multi-frame multiplier for the time-frequency lattice $\Z^d\times (D^{\tau})^{-1}\Z^d$. First, if $G(\x)$ is a Parseval Gabor multi-frame generator for $L^2(\R^d)$ and $H(\x)=M(\x)G(\x)$,
then  $H(\x)$ satisfies equation (\ref{cond12}), that is:
\begin{eqnarray*}
 d_0&=&  \sum_{\n\in\Z^d}\langle M(\x-\n)G(\x-\n), M(\x-\n)G(\x-\n)\rangle \\
&=& \sum_{\n\in\Z^d}\langle G(\x-\n), M^{*}(\x-\n)M(\x-\n)G(\x-\n)\rangle.
\end{eqnarray*}
Combining the above with
$$
d_0 = \sum_{\n\in\Z^d}\langle G(\x-\n), G(\x-\n)\rangle,$$
we obtain
\begin{eqnarray}\label{eq2.1}
\sum_{\n\in\Z^d}\langle G(\x-\n),(I-M^*(\x-\n)M(\x-\n))G(\x-\n)\rangle = 0.
\end{eqnarray}

\medskip
Now let  $\{\xi_1,\xi_2,...,\xi_\gamma\}$ be any orthonormal basis for $\C^\gamma$, and define
$$
G(\x)=\sum_{1\le j\le \gamma}\sqrt{d_0}\chi_{E_j}(\x)\xxi_j,
$$
where $E_1$, $E_2$, ..., $E_\gamma$ are as defined in Lemma \ref{Lem3}. For any $\x\in \O$, $\x\in E_j$ for some $j$, and equations (\ref{cond12}) and (\ref{cond22}) hold trivially for $G(\x)$ so it is a Parseval Gabor multi-frame generator for $L^2(\R^d)$ by Remark \ref{E_j_remark}. Furthermore, we have $G(\x-\n)=\0$ for any $\n\in \Z^d\setminus \{\0\}$ and $G(\x)=\xi_j$. Thus (\ref{eq2.1}) becomes
$$
d_0\langle \xi_j,(I-M^*(\x)M(\x))\xi_j\rangle = 0.
$$
Since $\xi_j$ can be any unit vector in $\C^\gamma$, this implies that $M^*(\x)M(\x)=I$ for any $x\in \O$ {\em a.e.}, hence for $x\in \R^d$ {\em a.e.} by Remark \ref{E_j_remark}.

\medskip\noindent
\textbf{Part 2.} Let $\k_0\in \KN$ be any vector. There are two cases to consider: $\k_0\not\in \L$ or $\k_0\in \L$.

\medskip\noindent
Case 1.  $\k_0\not\in \L$. Let $C$ be a small rectangle with the property described in and guaranteed by Lemma \ref{choice_lemma2}. For {a.e.} $x_0\in \R^d$, we can perform a translation so that $\x_0\in C^\p$ where $C^\p$ is the translation of $C$. Let $\O^t$ be  the  corresponding translation of $\O$, then  the previous discussions apply  to $\O^t$ and $C^\p$. Thus WLOG we can assume that $\x_0\in C$. We can choose $E_1$, ..., $E_\gamma$ such that $C\subset E_1$ and $-\k_0+C\subset E_2$ by Lemma \ref{choice_lemma2}. Define
$$
G(\x)=\sum_{1\le q\le \gamma}\sqrt{d_0}\chi_{E_q}(\x)\xi_q,
$$
where $\{\xi_1,\xi_2,...,\xi_\gamma\}$ is any orthonormal basis for $\R^\gamma$. (\ref{cond12}) and (\ref{cond22}) hold trivially, hence $G(\x)$ is a Parseval Gabor multi-frame generator for $L^2(\R^d)$. For $M(\x)G(\x)$ at $\x_0$ and $\k_0$, equation  (\ref{cond22}) contains only one term (since $\x_0-\k_0\in E_2$):
\begin{equation}\label{Meq}
\langle M(\x_0)G(\x_0),M(\x_0-\k_0)G(\x_0-\k_0)\rangle=
d_0\langle M(\x_0)\xi_1,M(\x_0-\k_0)\xi_2\rangle = 0.
\end{equation}
That is, $\xi_1^\tau M^*(\x_0)M(\x_0-\k_0)\xi_2=0$. Since $\{\xi_1,\xi_2,...,\xi_\gamma\}$ is arbitrary, we can replace $\xi_1$ and $\xi_2$ by $\e_i$ and $\e_j$ for any distinct $i$, $j$ between $1$ and $\gamma$ where $\{\e_1,\e_2,...,\e_d\}$ is the standard basis for $\R^d$. That is, $\e_i^\tau M^*(\x_0)M(\x_0-\k_0)\e_j=0$ for any $i\not=j$.
This  implies that $M^*(\x_0)M(\x_0-\k_0)$ is a diagonal matrix. On the other hand, if we replace $\xi_1$ and $\xi_2$ by $(\e_i+\e_j)/\sqrt{2}$ and $(\e_i-\e_j)/\sqrt{2}$ respectively, then equation (\ref{Meq}) leads to
 $\e_i^\tau M^*(\x_0)M(\x_0-\k_0)\e_i=\e_j^\tau M^*(\x_0)M(\x_0-\k_0)\e_j$ for any $i$ and $j$, proving that $M^*(\x_0)M(\x_0-\k_0)=\lambda_{\k_0}(\x_0)I$ with $\lambda_{\k_0}(\x_0)$ being  a unimodular scalar function.

\medskip\noindent
Case 2. $\k_0\in \L$. Let $C$ and $\k_{i_0}$ be as given in Lemma \ref{choice_lemma2} (2) so that $C\subset E_1$ and $-\k_{i_0}+C\subset E_2$. For any orthonormal basis $\{\xi_1,\xi_2,...,\xi_\gamma\}$ for $\R^\gamma$, define
\begin{eqnarray*}
G(\x)&=&\sqrt{d_0/2}\left(\chi_{E_1}(\x)\xi_1+\chi_{-\k_0+E_1}(\x)\xi_2+
\chi_{E_2}(\x)\xi_1-\chi_{-\k_0+E_2}(\x)\xi_2\right)\\
&+&\sum_{3\le q\le \gamma}\sqrt{d_0}\chi_{E_q}(\x)\xi_q.
\end{eqnarray*}
Equations (\ref{cond12}) and (\ref{cond22}) hold trivially for any $\x\in E_j$, $j\ge 3$, and (\ref{cond12}) also holds trivially for $\x\in E_1\cup E_2$. For $\x\in E_1$, and any $\k\in \KN$, (\ref{cond22}) contains only two non-trivial terms corresponding to $\l=\0$ and $\l=\k_0$, that is:
\begin{eqnarray}
&&\sum_{\l\in \L}\langle G(\x-\l), G(\x-\l-\k)\rangle\nonumber\\
&=&
\langle G(\x), G(\x-\k)\rangle+\langle G(\x-\k_0), G(\x-\k_0-\k)\rangle.\label{eq43}
\end{eqnarray}
If $\x-\k\not\in E_2$, then $\x-\k_0-\k\not\in -\k_0+E_2$ and both terms in (\ref{eq43}) equal to zero. If $\x-\k\in E_2$, then $\x-\k_0-\k\in -\k_0+E_2$ and (\ref{eq43}) becomes:
\begin{eqnarray*}
&&\langle G(\x), G(\x-\k)\rangle+\langle G(\x-\k_0), G(\x-\k_0-\k)\rangle\\
&=&
(d_0/2)\left(\langle \xi_1, \xi_1\rangle+\langle \xi_2, -\xi_2\rangle\right)=0.
\end{eqnarray*}
Thus (\ref{cond22}) holds for any $\x\in E_1$ and any $\k\in \KN$. Similarly, (\ref{cond22}) holds for any $\x\in E_2$ and any $\k\in \KN$. This proves that $G(\x)$ is a Parseval Gabor multi-frame generator for $L^2(\R^d)$. Now consider $M(\x)G(\x)$ at $\x_0$ and $\k_0\in \KN$. Equation  (\ref{cond22}) contains only one nontrivial term corresponding to $\l=\0$ (since the other possible non-trivial term corresponds to $\l=\k_0$ but $G(\x_0-\k_0-\k_0)=\0$):
$$
\langle M(\x_0)G(\x_0),M(\x_0-\k_0)G(\x_0-\k_0)\rangle=
(d_0/2)\langle M(\x_0)\xi_1,M(\x_0-\k_0)\xi_2\rangle = 0.
$$
That is, $\xi_1^\tau M^*(\x_0)M(\x_0-\k_0)\xi_2=0$. Since $\xi_1$ and $\xi_2$ are arbitrary, repeating the argument used in Case 1 leads to $M^*(\x_0)M(\x_0-\k_0)=\lambda_{\k_0}(\x_0)I$ with $\lambda_{\k_0}(\x_0)$ being  a unimodular scalar function.

\medskip\noindent
\textbf{Part 3.} Continue the discussion from Part 2 under the same setting and consider the two different cases.

\medskip\noindent
Case 1. $\k_0\not\in \L$. Recall that we have $C\subset E_1$ and $-\k_0+C\subset E_2$. For any given $\l_0\in \LN$, define
\begin{eqnarray*}
G(\x)&=&\sqrt{d_0/2}\left(\chi_{E_1}(\x)\e_1+\chi_{-\l_0+E_1}(\x)\e_2+
\chi_{E_2}(\x)\e_1-\chi_{-\l_0+E_2}(\x)\e_2\right)\\
&+&\sum_{3\le q\le \gamma}\sqrt{d_0}\chi_{E_q}(\x)\e_q.
\end{eqnarray*}
Again, equations (\ref{cond12}) and (\ref{cond22}) hold trivially for any $\x\in E_j$, $j\ge 3$, and (\ref{cond12}) also holds trivially for $\x\in E_1\cup E_2$. For any $\x\in E_1$, and any $\k\in \KN$, (\ref{cond22}) contains only two non-trivial terms corresponding to $\l=\0$ and $\l=\l_0$, that is:
\begin{eqnarray}
&&\sum_{\l\in \L}\langle G(\x-\l), G(\x-\l-\k)\rangle\nonumber\\
&=&
\langle G(\x), G(\x-\k)\rangle+\langle G(\x-\l_0), G(\x-\l_0-\k)\rangle.\label{eq44}
\end{eqnarray}
If $\x-\k\not\in E_2$, then $\x-\l_0-\k\not\in -\l_0+E_2$ and both terms in (\ref{eq44}) equal to zero.
If $\x-\k\in E_2$, then $\x-\l_0-\k\in -\l_0+E_2$ and (\ref{eq43}) becomes:
\begin{eqnarray*}
&&\langle G(\x), G(\x-\k)\rangle+\langle G(\x-\l_0), G(\x-\l_0-\k)\rangle\\
&=&
(d_0/2)\left(\langle \e_1, \e_1\rangle+\langle \e_2, -\e_2\rangle\right)=0.
\end{eqnarray*}
Thus (\ref{cond22}) holds for any $\x\in E_1$ and any $\k\in \KN$. Similarly, (\ref{cond22}) holds for any $\x\in E_2$ and any $\k\in \KN$. This proves that $G(\x)$ is a Parseval Gabor multi-frame generator for $L^2(\R^d)$. Now consider $M(\x)G(\x)$ at $\x_0$ and $\k_0\in \KN$. Equation  (\ref{cond22}) contains two nontrivial terms corresponding to $\l=\0$ and $\l=\l_0$ (since $\x_0-\k_0\in E_2$ and $\x_0-\l_0-\k_0\in -\l_0+E_2$):
\begin{eqnarray*}
&&\langle M(\x_0)G(\x_0),M(\x_0-\k_0)G(\x_0-\k_0)\rangle\\
&+&\langle M(\x_0-\l_0)G(\x_0-\l_0),M(\x_0-\l_0-\k_0)G(\x_0-\l_0-\k_0)\rangle\\
&=&
(d_0/2)\left(\langle M(\x_0)\e_1,M(\x_0-\k_0)\e_1\rangle -
\langle M(\x_0-\l_0)\e_2,M(\x_0-\l_0-\k_0)\e_2\rangle\right)\\
&=& 0.
\end{eqnarray*}
This implies that $\e_1^\tau M^*(\x_0)M(\x_0-\k_0)\e_1=\e_2^\tau M^*(\x_0-\l_0)M(\x_0-\l_0-\k_0)\e_2$. Since $M^*(\x_0)M(\x_0-\k_0)$ and $M^*(\x_0-\l_0)M(\x_0-\l_0-\k_0)$ are both scalar multiples of the identity matrix $I_{d\times d}$, this means that $M^*(\x_0)M(\x_0-\k_0)=M^*(\x_0-\l_0)M(\x_0-\l_0-\k_0)$ as desired.

\medskip\noindent
Case 2. $\k_0\in \L$. Recall that we have $C\subset E_1$ and $-\k_{i_0}+C\subset E_2$ for some $\k_{i_0}$ with $1\le i_0\le d$. For any given $\l_0\in \LN$, we need to consider several different cases.

\medskip\noindent
Subcase 1. $\l_0\not=\pm \k_0$. In this case we define
\begin{eqnarray*}
G(\x)&=&\sqrt{d_0}/2\left(\chi_{E_1}(\x)\e_1+\chi_{-\k_0+E_1}(\x)\e_1+
\chi_{-\l_0+E_1}(\x)\e_2-\chi_{-\l_0-\k_0+E_1}(\x)\e_2\right)\\
&+&\sqrt{d_0}/2\left(\chi_{E_2}(\x)\e_1-\chi_{-\k_0+E_2}(\x)\e_1+
\chi_{-\l_0+E_2}(\x)\e_2+\chi_{-\l_0-\k_0+E_2}(\x)\e_2\right)\\
&+&\sum_{3\le q\le \gamma}\sqrt{d_0}\chi_{E_q}(\x)\e_q.
\end{eqnarray*}
Equations (\ref{cond12}) and (\ref{cond22}) hold trivially for any $\x\in E_j$, $j\ge 3$, and (\ref{cond12}) also holds trivially for $\x\in E_1\cup E_2$. For any $\x\in E_1$, and any $\k\in \KN$, (\ref{cond22}) contains only four non-trivial terms corresponding to $\l=\0$, $\l=\k_0$, $\l=\l_0$ and $\l=\l_0+\k_0$, that is:
\begin{eqnarray}
&&\sum_{\l\in \L}\langle G(\x-\l), G(\x-\l-\k)\rangle\nonumber\\
&=&
\langle G(\x), G(\x-\k)\rangle+\langle G(\x-\k_0), G(\x-\k_0-\k)\rangle\nonumber\\
&+&
\langle G(\x-\l_0), G(\x-\l_0-\k)\rangle +\langle G(\x-\l_0-\k_0), G(\x-\l_0-\k_0-\k)\rangle\nonumber\\
&=&
\sqrt{d_0}/2\langle \e_1, G(\x-\k)\rangle+\sqrt{d_0}/2\langle \e_1, G(\x-\k_0-\k)\rangle\nonumber\\
&+&
\sqrt{d_0}/2\langle \e_2, G(\x-\l_0-\k)\rangle +\sqrt{d_0}/2\langle -\e_2, G(\x-\l_0-\k_0-\k)\rangle.\label{eq45}
\end{eqnarray}
If $\x-\k\not\in E_2$ and $\x-\k\not\in -\k_0+E_2$, then $\x-\l_0-\k\not\in -\l_0+E_2$ and $\x-\l_0-\k_0-\k\not\in -\l_0-\k_0+E_2$ and each term in (\ref{eq45}) equals zero.
If $\x-\k\in E_2$, then $\x-\k_0-\k\in -\k_0+E_2$, $\x-\l_0-\k\in -\l_0+E_2$ and $\x-\l_0-\k_0-\k\in -\l_0-\k_0+E_2$. Thus (\ref{eq45}) becomes:
$$
(d_0/4)\left(\langle \e_1, \e_1\rangle+\langle \e_1, -\e_1\rangle
+\langle \e_2, \e_2\rangle+\langle -\e_2, \e_2\rangle\right)=0.
$$
If $\x-\k\in -\k_0+E_2$, then $\x-\k_0-\k\in -2\k_0+E_2$, $\x-\l_0-\k\in -\l_0-\k_0+E_2$ and $\x-\l_0-\k_0-\k\in -\l_0-2\k_0+E_2$. Since $-2\k_0+E_2$ is disjoint from $E_2$ and $-\k_0+E_2$, $G(\x-\k_0-\k)\not=(\sqrt{d_0}/2)\e_1$. Similarly, $-\l_0-2\k_0+E_2$ is disjoint from $-\l_0+E_2$ and $-\l_0-\k_0+E_2$ hence $G(\x-\l_0-\k_0-\k)\not=(\sqrt{d_0}/2)\e_2$. It follows that
(\ref{eq45}) contains only two nontrivial terms corresponding to $G(\x-\k)$ and $G(\x-\l_0-\k)$, which   becomes:
$$
(d_0/4)\left(\langle \e_1, -\e_1\rangle+\langle \e_2, \e_2\rangle\right)=0.
$$
The case $\x\in E_2$ can  be similarly verified. Thus $G(\x)$ is a Parseval Gabor multi-frame generator for $L^2(\R^d)$. Substituting $G(\x)$ in (\ref{cond22}) by $M(\x)G(\x)$ with $\x=\x_0$ and $\k=\k_0$ then yields (keep in mind that $G(\x_0-2\k_0)=\0$ and $G(\x_0-\l_0-2\k_0)=\0$):
\begin{eqnarray*}
&&\langle M(\x_0)G(\x_0),M(\x_0-\k_0)G(\x_0-\k_0)\rangle\\
&+&\langle M(\x_0-\l_0)G(\x_0-\l_0),M(\x_0-\l_0-\k_0)G(\x_0-\l_0-\k_0)\rangle\\
&=&
(d_0/4)\left(\langle M(\x_0)\e_1,M(\x_0-\k_0)\e_1\rangle -
\langle M(\x_0-\l_0)\e_2,M(\x_0-\l_0-\k_0)\e_2\rangle\right)\\
&=& 0.
\end{eqnarray*}
This implies that $\e_1^\tau M^*(\x_0)M(\x_0-\k_0)\e_1=\e_2^\tau M^*(\x_0-\l_0)M(\x_0-\l_0-\k_0)\e_2$. This means that $M^*(\x_0)M(\x_0-\k_0)=M^*(\x_0-\l_0)M(\x_0-\l_0-\k_0)$ as desired.

\medskip\noindent
Subcase 2. $\l_0=\k_0$. In this case we define
\begin{eqnarray*}
G(\x)&=&\sqrt{d_0}/2\left(\chi_{E_1}(\x)\e_1+\chi_{-\l_0+E_1}(\x)(-\e_1+\e_2)+
\chi_{-2\l_0+E_1}(\x)\e_2\right)\\
&+&\sqrt{d_0}/2\left(\chi_{E_2}(\x)\e_1+\chi_{-\l_0+E_2}(\x)(\e_1+\e_2)-
\chi_{-2\l_0+E_2}(\x)\e_2\right)\\
&+&\sum_{3\le q\le \gamma}\sqrt{d_0}\chi_{E_q}(\x)\e_q.
\end{eqnarray*}
We leave it to our reader to verify that equations (\ref{cond12}) and (\ref{cond22}) hold for any $\x\in \O$ (hence for any $\x\in \R^d$) and for any $\k\in \KN$, that is, $G(\x)$ is a Parseval Gabor multi-frame generator for $L^2(\R^d)$. Substituting $G(\x)$ in (\ref{cond22}) by $M(\x)G(\x)$ with $\x=\x_0$ and $\k=\k_0$ then yields
\begin{eqnarray*}
&&\langle M(\x_0)G(\x_0),M(\x_0-\k_0)G(\x_0-\k_0)\rangle\\
&+&\langle M(\x_0-\l_0)G(\x_0-\l_0),M(\x_0-\l_0-\k_0)G(\x_0-\l_0-\k_0)\rangle\\
&=&
(d_0/4)\left(\langle M(\x_0)\e_1,M(\x_0-\k_0)(-\e_1+\e_2)\rangle\right)\\
&+&
(d_0/4)\left(\langle M(\x_0-\l_0)(-\e_1+\e_2),M(\x_0-\l_0-\k_0)\e_2\rangle\right)\\
&=& 0.
\end{eqnarray*}
That is,
\begin{eqnarray*}
&&\e_1^\tau M^*(\x_0)M(\x_0-\k_0)(-\e_1+\e_2)\\
&=&
-(-\e_1+\e_2)^\tau M^*(\x_0-\l_0)M(\x_0-\l_0-\k_0)\e_2.
\end{eqnarray*}

 Since $M^*(\x_0)M(\x_0-\k_0)$ and $M^*(\x_0-\l_0)M(\x_0-\l_0-\k_0)$ are both scalar multiples of the identity matrix $I_{d\times d}$, the above simplifies to
 $$
 \e_1^\tau M^*(\x_0)M(\x_0-\k_0)\e_1=\e_2^\tau M^*(\x_0-\l_0)M(\x_0-\l_0-\k_0)\e_2.
 $$
 This means that $M^*(\x_0)M(\x_0-\k_0)=M^*(\x_0-\l_0)M(\x_0-\l_0-\k_0)$ as desired.

\medskip\noindent
Subcase 3. $\l_0=-\k_0$. In this case we define $G(\x)$ by
\begin{eqnarray*}
G(\x)&=&\sqrt{d_0}/2\left(\chi_{\l_0+E_1}(\x)\e_1+\chi_{E_1}(\x)(-\e_1+\e_2)+
\chi_{-\l_0+E_1}(\x)\e_2\right)\\
&+&\sqrt{d_0}/2\left(\chi_{\l_0+E_2}(\x)\e_1+\chi_{E_2}(\x)(\e_1+\e_2)-
\chi_{-\l_0+E_2}(\x)\e_2\right)\\
&+&\sum_{3\le q\le \gamma}\sqrt{d_0}\chi_{E_q}(\x)\e_q.
\end{eqnarray*}
The rest of the proof is similar to Subcase 2 and is left to the reader.

\medskip
Since $\k_0\in \K$, $\l_0\in \L$ are arbitrary and $\x_0$ is any point in $\R^d$ (in the {\em a.e.} sense), we have shown that $M^*(\x)M(\x-\k)$ is $\L$-periodic for any $\x\in \R^d$ {\em a.e.} and any $\k\in \KN$. This concludes the proof of Theorem \ref{main-thm} and also this paper.

\end{document}